\documentclass[12pt,twoside]{report}
\usepackage{svcon2e}
\usepackage{amsthm,amssymb}
\usepackage{amsbsy,amsfonts,amsmath,verbatim}

\catcode`@=12
\makeatletter

\@addtoreset{equation}{section}

\theoremstyle{plain}

\newtheorem{theorem}{Theorem}[section]

\newtheorem{@definition}{\bf Definition}[section]

\newtheorem{@example}{\bf Example}[section]
\newenvironment{example}{\begin{@example}\rm}{\end{@example}}

\newtheorem{@remark}{\bf Remark}[section]

\newcommand{\R}{R}
\newcommand{\Rt}{\tilde{R}}
\newcommand{\RC}{\mathrm{RC}}
\newcommand{\RH}{\mathcal{R}}
\newcommand{\St}{\tilde{S}}
\newcommand{\Tab}{\mathrm{Tab}}
\newcommand{\Y}{\mathcal{Y}}
\newcommand{\Z}{\mathbb{Z}}
\newcommand{\bins}[2]{{\textstyle\genfrac{(}{)}{0pt}{}{#1}{#2}}}
\newcommand{\bin}[2]{\genfrac{(}{)}{0pt}{}{#1}{#2}}
\newcommand{\card}{\mathrm{card}}
\newcommand{\col}{\mathrm{col}}
\newcommand{\cospin}{\mathrm{cospin}}
\newcommand{\ct}{\tilde{c}}
\newcommand{\diag}{\mathrm{diag}}
\newcommand{\height}{h}
\newcommand{\inv}{\mathrm{inv}}
\newcommand{\la}{\lambda}

\newcommand{\lecell}{\preceq}
\newcommand{\lesscell}{\prec}
\newcommand{\maxspin}{\mathrm{maxspin}}
\newcommand{\mud}{\mu^\bullet}
\newcommand{\nud}{\nu^\bullet}
\newcommand{\pos}{\mathrm{pos}}
\newcommand{\ribcov}{\lessdot_L}
\newcommand{\rible}{\le_L}
\newcommand{\row}{\mathrm{row}}
\newcommand{\spa}{\phantom{1}}
\newcommand{\spin}{\mathrm{spin}}
\newcommand{\st}{\mathrm{st}}
\newcommand{\up}[1]{#1\!\uparrow}
\newcommand{\qbins}[2]{{\textstyle\genfrac{[}{]}{0pt}{}{#1}{#2}}}
\newcommand{\qbin}[2]{\genfrac{[}{]}{0pt}{}{#1}{#2}}
\newcommand{\quot}{\mathrm{quot}}

\begin{document}
\pagenumbering{arabic}

\chapter{$q$-Supernomial Coefficients: From Riggings to Ribbons}
\chapterauthors{Anne Schilling}

\vspace{0.4cm}
\begin{center}
Dedicated to Barry M. McCoy on the occasion of his sixtieth birthday
\end{center}

\begin{abstract}
$q$-Supernomial coefficients are generalizations of the $q$-binomial
coefficients. They can be defined as the coefficients of the 
Hall--Littlewood symmetric function in a product of the complete 
symmetric functions or the elementary symmetric functions. 
Hatayama et al. 
give an explicit expression for these $q$-supernomial coefficients. 
A combinatorial expression as the generating function of ribbon
tableaux with (co)spin statistic follows from the work of
Lascoux, Leclerc and Thibon. 
In this paper we interpret the formulas by Hatayama et al. in terms of 
rigged configurations and provide an explicit statistic preserving
bijection between rigged configurations and ribbon tableaux thereby
establishing a new direct link between these combinatorial objects.
\end{abstract}

\section{Introduction}\label{sec:intro}

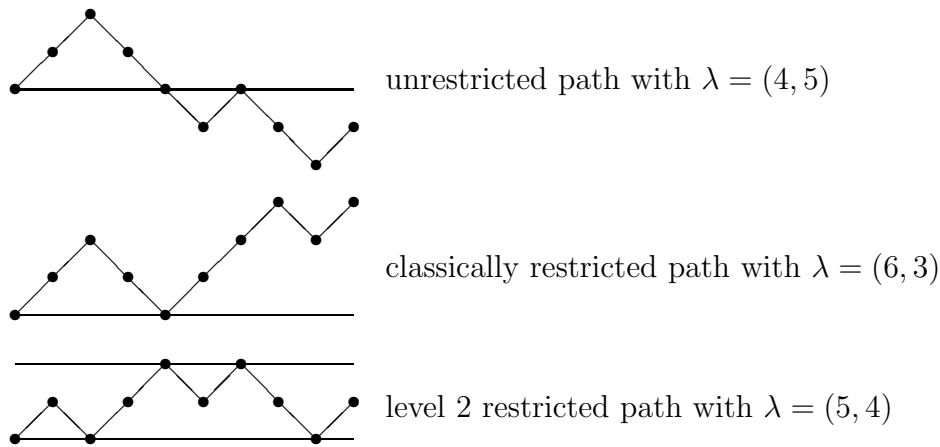
\begin{figure}
\unitlength0.5cm
\begin{tabular}{ll}
\raisebox{-1cm}{
\begin{picture}(9,4)
\put(0,2){\circle*{0.3}}
\put(0,2){\line(1,1){2}}
\put(1,3){\circle*{0.3}}
\put(2,4){\circle*{0.3}}
\put(2,4){\line(1,-1){3}}
\put(3,3){\circle*{0.3}}
\put(4,2){\circle*{0.3}}
\put(5,1){\circle*{0.3}}
\put(5,1){\line(1,1){1}}
\put(6,2){\circle*{0.3}}
\put(6,2){\line(1,-1){2}}
\put(7,1){\circle*{0.3}}
\put(8,0){\circle*{0.3}}
\put(8,0){\line(1,1){1}}
\put(9,1){\circle*{0.3}}
\put(0,2){\line(1,0){9}}
\end{picture}}
& unrestricted path with $\la=(4,5)$\\[5mm]
\raisebox{-0.5cm}{
\begin{picture}(9,4)
\put(0,0){\circle*{0.3}}
\put(0,0){\line(1,1){2}}
\put(1,1){\circle*{0.3}}
\put(2,2){\circle*{0.3}}
\put(2,2){\line(1,-1){2}}
\put(3,1){\circle*{0.3}}
\put(4,0){\circle*{0.3}}
\put(4,0){\line(1,1){3}}
\put(5,1){\circle*{0.3}}
\put(6,2){\circle*{0.3}}
\put(7,3){\circle*{0.3}}
\put(7,3){\line(1,-1){1}}
\put(8,2){\circle*{0.3}}
\put(8,2){\line(1,1){1}}
\put(9,3){\circle*{0.3}}
\put(0,0){\line(1,0){9}}
\end{picture}}
& classically restricted path with $\la=(6,3)$\\[10mm]
\raisebox{-0.3cm}{
\begin{picture}(9,2)
\put(0,0){\circle*{0.3}}
\put(1,1){\circle*{0.3}}
\put(0,0){\line(1,1){1}}
\put(1,1){\line(1,-1){1}}
\put(2,0){\circle*{0.3}}
\put(2,0){\line(1,1){2}}
\put(3,1){\circle*{0.3}}
\put(4,2){\circle*{0.3}}
\put(4,2){\line(1,-1){1}}
\put(5,1){\circle*{0.3}}
\put(5,1){\line(1,1){1}}
\put(6,2){\circle*{0.3}}
\put(6,2){\line(1,-1){2}}
\put(7,1){\circle*{0.3}}
\put(8,0){\circle*{0.3}}
\put(8,0){\line(1,1){1}}
\put(9,1){\circle*{0.3}}
\put(0,0){\line(1,0){9}}
\put(0,2){\line(1,0){9}}
\end{picture}}
& level 2 restricted path with $\la=(5,4)$
\end{tabular}
\vspace{0.5cm}
\caption{\label{fig:paths}Examples of unrestricted, classically 
restricted and level-restricted paths.}
\end{figure}

Lattice paths play an important r\^{o}le in combinatorics and
exactly solvable lattice models. One distinguishes three types
of paths: unrestricted, classically restricted and level-restricted
paths. Amongst the easiest examples are Dyck paths consisting
of up- and down-steps. For fixed $\la=(\la_1,\la_2)$, the set
of unrestricted paths contains all paths with $\la_1$ up-steps and
$\la_2$ down-steps. A path is classically restricted if the number
of up-steps is greater or equal to the the number of down-steps in the 
first $k$ steps for all $1\le k\le \la_1+\la_2$. A path is restricted
of level $\ell$ if it is classically restricted and in addition the 
difference between the number of up- and down-steps in the first $k$ 
steps does not exceed $\ell$ for all $1\le k\le \la_1+\la_2$.
Examples of all three types of paths can be found in Figure \ref{fig:paths}.

The number of unrestricted paths with $\la_1$ up-steps and
$\la_2$ down-steps is given by the binomial coefficient 
$\bins{\la_1+\la_2}{\la_1,\la_2}=(\la_1+\la_2)!/\la_1! \la_2!$
which is the expansion coefficient of
\begin{equation}\label{eq:bin}
(x_1+x_2)^L = \sum_{\substack{\la_1,\la_2\ge 0\\ \la_1+\la_2=L}}
 x_1^{\la_1} x_2^{\la_2} \bin{L}{\la_1,\la_2}.
\end{equation}

More generally, the steps of paths associated with the Lie algebra 
$A_{n-1}$ are Young tableaux over the alphabet 
$\{1,2,\ldots,n\}$ \cite{JMO:1988}.
A Young tableau is a filling of a partition shape 
$\tau=(\tau_1\ge \tau_2 \ge \cdots \ge \tau_n\ge 0)$ which is weakly
increasing along rows and strictly increasing along columns.
The content of a tableau $\la=(\la_1,\ldots,\la_n)$
records the number occurrences of the various letters in the
tableau, i.e., $\la_i$ specifies the numbers of $i$'s in the tableau.
In this language, Dyck paths are associated with the algebra
$A_1$, and the up- and down-steps correspond to single-box Young tableaux 
filled with either 1 or 2. 

The number of paths with single-row steps of widths
$\mu_1,\mu_2,\ldots,\mu_L$ and total content 
$\la=(\la_1,\ldots,\la_n)$, denoted by $S_{\la\mu}$, is given 
by the coefficient of $x^\la:=x_1^{\la_1}x_2^{\la_2}\cdots x_n^{\la_n}$
in the expansion of $h_{\mu_1}\cdots h_{\mu_L}$
where $h_r$ is the complete symmetric polynomial
of degree $r$ in $n$ variables \cite[Chap. I.2]{Mac:1995}.
In analogy to \eqref{eq:bin} we have
\begin{equation*}
h_{\mu_1}\cdots h_{\mu_L}=\sum_{\la} x^\la S_{\la\mu}.
\end{equation*}
Since the $S_{\la\mu}$ generalize the binomial and multinomial 
coefficients, they were coined supernomial coefficients 
(or more precisely completely symmetric supernomial coefficients)
in \cite{SW:1998,SW:1999}.

Similarly, the number of paths with single-column steps
with heights $\mu_1,\mu_2,\ldots,\mu_L$ and total content 
$\la=(\la_1,\ldots,\la_n)$, denoted by $S'_{\la\mu}$, is given 
by the coefficient of $x^\la:=x_1^{\la_1}x_2^{\la_2}\cdots x_n^{\la_n}$
in the expansion of $e_{\mu_1}\cdots e_{\mu_L}$
where $e_r$ is the elementary symmetric polynomial
of degree $r$ in $n$ variables \cite[Chap. I.2]{Mac:1995}.
That is
\begin{equation*}
e_{\mu_1}\cdots e_{\mu_L}=\sum_{\la} x^\la S'_{\la\mu}.
\end{equation*}
$S'_{\la\mu}$ is called the completely antisymmetric supernomial
coefficient.

The one-dimensional configuration sums of exactly solvable
lattice models, which are necessary for the calculation of order parameters
of these models, require $q$-analogues of the supernomial coefficients.
For example, the $q$-analogue of the binomial coefficient is the
$q$-binomial coefficient
\begin{equation}\label{eq:qbin}
\qbin{\la_1+\la_2}{\la_1,\la_2} = \begin{cases}
 \frac{(q)_{\la_1+\la_2}} {(q)_{\la_1}(q)_{\la_2}}
 &\text{if $\la_1,\la_2$ are nonnegative integers,}\\
 0 &\text{otherwise,}
\end{cases}
\end{equation}
where $(q)_m=\prod_{i=1}^m (1-q^i)$. 

The $q$-analogues of $S_{\la\mu}$ and $S'_{\la\mu}$ are the
coefficients of the Hall--Littlewood function $P_\mu(x;q)$
\cite[Chap. III.2]{Mac:1995} in the expansion of 
$h_{\la_1}(x)\cdots h_{\la_n}(x)$ and $e_{\la_1}(x)\cdots e_{\la_n}(x)$,
respectively. Here $h_r(x)$ and $e_r(x)$ are the complete and
elementary symmetric functions in infinitely many variables.
The $q$-supernomial coefficients can be expressed in terms of the 
Kostka polynomials $K_{\la\mu}(q)$ which are the entries of the 
transition matrix between the Schur function $s_\la(x)$ 
\cite[Chaps. I.3]{Mac:1995} and the Hall--Littlewood function $P_\mu(x;q)$
\begin{equation}\label{eq:HL}
s_\la(x) = \sum_\mu K_{\la\mu}(q) P_\mu(x;q).
\end{equation}
Combining \eqref{eq:HL} with the expansions
$h_{\la_1}(x)\cdots h_{\la_n}(x)=\sum_\eta K_{\eta\la} s_\eta(x)$
and $e_{\la_1}(x)\cdots e_{\la_n}(x)=\sum_\eta K_{\eta^t\la} s_\eta(x)$
yields
\begin{align}
S_{\la\mu}(q)=\sum_\eta K_{\eta\la} K_{\eta\mu}(q)\\
\intertext{and}
S'_{\la\mu}(q)=\sum_\eta K_{\eta^t\la} K_{\eta\mu}(q).
\end{align}
Here $K_{\la\mu}:=K_{\la\mu}(1)$ are the Kostka numbers and
$\eta^t$ is the transpose partition of $\eta$ obtained by
interchanging the r\^{o}le of rows and columns.

The $q$-analogues of the trinomial coefficients, which correspond
to the case $\la=(\la_1,\la_2)$ and $\mu=(2^L)$, were extensively
studied by Andrews and Baxter \cite{Andrews:1994,AB:1987}.
For $\la=(\la_1,\la_2)$ and general $\mu$ explicit formulas for 
$S_{\la\mu}(q)$ were given in \cite{SW:1998} where they were also 
used to prove generalizations of Rogers--Ramanujan-type identities.
Recently, an alternative definition as characters of coinvariants for one 
dimensional vertex operator algebras was given in \cite{FLT:2001}.
Hatayama et al. \cite{HKKOTY:1999,Kirillov:2000} provide an explicit expression
for $S_{\la\mu}(q)$ and $S'_{\la\mu}(q)$ as given in equations 
\eqref{eq:sym-super} and \eqref{eq:anti-super} below. 
The motivation for these formulas comes from exactly solvable models,
where $S_{\la\mu}(q)$ and $S'_{\la\mu}(q)$ can be interpreted as the 
generating function of unrestricted paths with an energy statistic which comes
from crystal base theory.

Lascoux, Leclerc and Thibon \cite{LLT:1997} introduced
$q$-supernomial coefficients as the generating function of
ribbon tableaux with (co)spin statistic. As was shown in
\cite{LT:2000} these $q$-analogues are related to parabolic
Kazhdan--Lusztig polynomials for affine symmetric groups.
There exists a bijection between $L$-ribbon tableaux and $L$-tuples
of ordinary tableaux \cite{SW:1985}.
Recently \cite{SSW:2000}, the (co)spin statistic was translated into
inversion statistic on tuples of tableaux under this bijection.

In this paper we show that the expressions \eqref{eq:sym-super} 
and \eqref{eq:anti-super} for the $q$-supernomial coefficients
given below have a combinatorial interpretation in terms of 
``rigged configurations'' similar to the ones introduced by Kerov, Kirillov 
and Reshetikhin \cite{KKR:1988,KR:1988,KR:1990} for generating functions 
of classically restricted paths.
In addition we give a statistic preserving bijection between these
rigged configurations and ribbon tableaux, thereby providing a new
direct link between these combinatorial objects.

The paper is organized as follows. In Section \ref{sec:RC} we 
state the explicit formulas for the $q$-supernomial coefficients
as given by Hatayama et al. \cite{HKKOTY:1999} and provide the combinatorial
description in terms of rigged configurations.
In Section \ref{sec:RT} we review the generating functions
of $L$-ribbon tableaux with (co)spin statistic and the analogous
interpretation in terms of inversion statistic on $L$-tuples of
tableaux. The bijection between ribbon tableaux and rigged 
configurations is presented in Section \ref{sec:bij}.

\section{Rigged configurations}
\label{sec:RC}

In this section we present the explicit expressions for the
$q$-supernomial coefficients $S_{\la\mu}(q)$ and $S'_{\la\mu}(q)$
of Hatayama et al. \cite{HKKOTY:1999} and provide a combinatorial 
interpretation of these formulas in terms of rigged configurations.

\subsection{Explicit expressions for the $q$-supernomials}

In the following it will be useful to consider a slight variant
of the $q$-supernomial coefficients $S_{\la\mu}(q)$ and
$S'_{\la\mu}(q)$ using a costatistic. Let $n(\mu)=\sum_{1\le i<j\le L}
\min(\mu_i,\mu_j)$ for $\mu=(\mu_1,\ldots,\mu_L)$. Then define
\begin{equation*}
\St_{\la\mu}(q)=q^{n(\mu)} S_{\la\mu}(q^{-1}) \qquad \text{and} \qquad 
\St'_{\la\mu}(q)=q^{n(\mu)} S'_{\la\mu}(q^{-1}).
\end{equation*}

Let $\la=(\la_1,\ldots,\la_n)$ and $\mu=(\mu_1,\ldots,\mu_L)$ be 
partitions. According to \cite{HKKOTY:1999,Kirillov:2000}, the 
completely symmetric $q$-supernomial is given by
\begin{equation}\label{eq:sym-super}
\St_{\la\mu}(q)=\sum_{\{\nu\}} q^{\Phi(\nu)} 
 \prod_{\substack{1\le a\le n-1\\ 1\le i\le \mu_1}}
 \qbin{\nu_i^{(a+1)}-\nu_{i+1}^{(a)}}{\nu_i^{(a)}-\nu_{i+1}^{(a)},
 \nu_i^{(a+1)}-\nu_i^{(a)}},
\end{equation}
where the sum $\sum_{\{\nu\}}$ is over all sequences of partitions
$\nu^{(1)},\ldots,\nu^{(n-1)}$ such that
\begin{equation}
\begin{split}
&\emptyset=\nu^{(0)}\subset \nu^{(1)} \subset \cdots \subset \nu^{(n-1)}
 \subset \nu^{(n)} = \mu^t\\
&|\nu^{(a)}|=\la_1+\cdots+\la_a \quad \text{for all $1\le a<n$}
\end{split}
\end{equation}
and
\begin{equation}
\Phi(\nu)=\sum_{\substack{1\le a <n\\ 1\le i\le \mu_1}}
 \nu_{i+1}^{(a)}(\nu_i^{(a+1)}-\nu_i^{(a)}).
\end{equation}

\begin{example}\label{ex:config sym}
Take $\la=\mu=(2,2,1)$. Then the allowed sequences $\nu$
and their contributions to the sum in \eqref{eq:sym-super} are given by
\begin{equation*}
\begin{array}{llllllll}
 \emptyset 
&\subset
&\raisebox{-0.2cm}{\begin{array}[b]{|c|c|} \hline \spa & \spa \\ 
  \hline \end{array}}
&\subset
&\raisebox{-0.2cm}{\begin{array}[b]{|c|c|} \hline \spa & \spa \\ 
  \hline \spa & \spa \\ \hline \end{array}}
&\subset
&\raisebox{-0.2cm}{\begin{array}[b]{|c|c|c|} \cline{1-2} \spa & \spa & 
  \multicolumn{1}{c}{} \\ \hline \spa & \spa & \spa \\ \hline \end{array}}
& \qquad q^2 \\[4mm]
 \emptyset 
&\subset
&\raisebox{-0.2cm}{\begin{array}[b]{|c|c|} \hline \spa & \spa \\ 
  \hline \end{array}}
&\subset
&\raisebox{-0.2cm}{\begin{array}[b]{|c|c|c|} 
  \cline{1-1} \spa & \multicolumn{2}{c}{} \\
  \hline \spa & \spa & \spa \\ \hline \end{array}} 
&\subset
&\raisebox{-0.2cm}{\begin{array}[b]{|c|c|c|} \cline{1-2} \spa & \spa & 
  \multicolumn{1}{c}{} \\ \hline \spa & \spa & \spa \\ \hline \end{array}}
& \qquad \phantom{q^0} \qbin{3}{2,1}\qbin{2}{1,1} \\[4mm]
 \emptyset 
&\subset
&\raisebox{-0.2cm}{\begin{array}[b]{|c|} \hline \spa \\ \hline \spa \\ 
  \hline \end{array}}
&\subset
&\raisebox{-0.2cm}{\begin{array}[b]{|c|c|c|} 
  \cline{1-1} \spa & \multicolumn{2}{c}{} \\
  \hline \spa & \spa & \spa \\ \hline \end{array}} 
&\subset
&\raisebox{-0.2cm}{\begin{array}[b]{|c|c|c|} \cline{1-2} \spa & \spa & 
  \multicolumn{1}{c}{} \\ \hline \spa & \spa & \spa \\ \hline \end{array}}
 & \qquad q^2 \qbin{2}{1,1} \\[4mm]
 \emptyset 
&\subset
&\raisebox{-0.2cm}{\begin{array}[b]{|c|} \hline \spa \\ \hline \spa \\ 
  \hline \end{array}}
&\subset
&\raisebox{-0.2cm}{\begin{array}[b]{|c|c|} \hline \spa & \spa \\ 
  \hline \spa & \spa \\ \hline \end{array}}
&\subset
&\raisebox{-0.2cm}{\begin{array}[b]{|c|c|c|} \cline{1-2} \spa & \spa & 
  \multicolumn{1}{c}{} \\ \hline \spa & \spa & \spa \\ \hline \end{array}}
 & \qquad q^3 \qbin{2}{1,1}
\end{array}
\end{equation*}
Hence $\St_{\la\mu}(q)=1+2q+4q^2+3q^3+q^4$.
\end{example}

Similarly, the explicit expression for the completely antisymmetric 
supernomial is given by \cite{HKKOTY:1999,Kirillov:2000}
\begin{equation}\label{eq:anti-super}
\St'_{\la\mu}(q)=\sum_{\{\nu\}} 
 \prod_{\substack{1\le a\le n-1\\ 1\le i\le \mu_1}}
 \qbin{\nu_i^{(a+1)}-\nu_{i+1}^{(a+1)}}{\nu_i^{(a)}-\nu_{i+1}^{(a+1)},
\nu_i^{(a+1)}-\nu_i^{(a)}},
\end{equation}
where the sum $\sum_{\{\nu\}}$ runs over all sequences of partitions
$\nu^{(1)},\ldots,\nu^{(n-1)}$ such that
\begin{equation}
\begin{split}
&\emptyset=\nu^{(0)}\subset \nu^{(1)} \subset \cdots \subset \nu^{(n-1)}
 \subset \nu^{(n)} = \mu^t\\
&\nu^{(a)}/\nu^{(a-1)} \quad \text{is a horizontal $\la_a$-strip.}
\end{split}
\end{equation}
Here $\nu^{(a)}/\nu^{(a-1)}$ is a skew shape obtained by considering
all boxes in $\nu^{(a)}$ not in $\nu^{(a-1)}$, and a horizontal $p$-strip 
is a skew shape with $p$ boxes such that every column contains at most
one box.

\begin{example}\label{ex:config anti}
For the antisymmetric case also take $\la=\mu=(2,2,1)$. Then the allowed 
sequences $\nu$ and their contributions to the sum in \eqref{eq:anti-super}
are given by
\begin{equation*}
\begin{array}{llllllll}
 \emptyset 
&\subset
&\raisebox{-0.2cm}{\begin{array}[b]{|c|c|} \hline \spa & \spa \\ 
  \hline \end{array}}
&\subset
&\raisebox{-0.2cm}{\begin{array}[b]{|c|c|} \hline \spa & \spa \\ 
  \hline \spa & \spa \\ \hline \end{array}}
&\subset
&\raisebox{-0.2cm}{\begin{array}[b]{|c|c|c|} \cline{1-2} \spa & \spa & 
  \multicolumn{1}{c}{} \\ \hline \spa & \spa & \spa \\ \hline \end{array}}
& \qquad 1 \\[4mm]
 \emptyset 
&\subset
&\raisebox{-0.2cm}{\begin{array}[b]{|c|c|} \hline \spa & \spa \\ 
  \hline \end{array}}
&\subset
&\raisebox{-0.2cm}{\begin{array}[b]{|c|c|c|} 
  \cline{1-1} \spa & \multicolumn{2}{c}{} \\
  \hline \spa & \spa & \spa \\ \hline \end{array}} 
&\subset
&\raisebox{-0.2cm}{\begin{array}[b]{|c|c|c|} \cline{1-2} \spa & \spa & 
  \multicolumn{1}{c}{} \\ \hline \spa & \spa & \spa \\ \hline \end{array}}
& \qquad \qbin{2}{1,1}^2
\end{array}
\end{equation*}
Hence $\St'_{\la\mu}(q)=2+2q+q^2$.
\end{example}

\subsection{Combinatorial interpretation in terms of rigged configurations}

The main tool for the combinatorial interpretation of \eqref{eq:sym-super} 
and \eqref{eq:anti-super} is the interpretation of the $q$-binomial coefficient
$\qbins{m+p}{m,p}$ as the generating function of partitions in a
box of width $m$ and height $p$ (see \cite[Theorem 3.1]{Andrews:1976}).
That is
\begin{equation*}
\qbin{m+p}{m,p} = \sum_{\la\subset (m^p)} q^{|\la|}
\end{equation*}
where $|\la|=\la_1+\la_2+\cdots$ for the partition $\la=(\la_1,\la_2,\ldots)$.

For the completely symmetric $q$-supernomials $\St_{\la\mu}(q)$, set
$m_i^{(a)}=\nu_i^{(a)}-\nu_{i+1}^{(a)}$ and
$p_i^{(a)}=\nu_i^{(a+1)}-\nu_i^{(a)}$. The $p_i^{(a)}$ are called
vacancy numbers. For fixed $\nu$, each term in the summand 
of \eqref{eq:sym-super} corresponds to a labeling or rigging of
$\nu$ in the following way. Label each column of $\nu^{(a)}$
with a quantum number $j$. If the column in $\nu^{(a)}$ has height
$i$, the quantum number $j$ has to be an integer satisfying
$0\le j \le p_i^{(a)}$. If $j=p_i^{(a)}$ it is called a singular
quantum number. Riggings which differ only by reordering of quantum 
numbers corresponding to columns of the same height in a partition are 
identified. Hence for each partition $\nu^{(a)}$
and column height $i$ we may view the riggings as a partition
$J_i^{(a)}$ in a box of width $m_i^{(a)}$ and height $p_i^{(a)}$.
A configuration together with an admissible rigging is called a rigged 
configuration.

\begin{example}
Continuing example \ref{ex:config sym},
\begin{equation*}
\begin{array}[b]{|c|c|l} \cline{1-2} 1 & 0 & 1 \\ \cline{1-2} \end{array} 
\qquad
\begin{array}[b]{|c|c|c|l} \cline{1-1} 1 & \multicolumn{3}{l}{1} \\
 \cline{1-3} \spa & 0 & 0 & 0 \\ \cline{1-3} \end{array}
\end{equation*}
is a valid rigging of the second sequence, where the vacancy number
$p_i^{(a)}$ is written on the right of row $i$ in $\nu^{(a)}$ and the 
quantum numbers are put in the highest box of each column.
\end{example}

Let us denote the set of all rigged configurations corresponding
to $\la$ and $\mu$ by $\RC(\la,\mu)$. To each $(\nu,J)\in \RC(\la,\mu)$
associate a statistic
\begin{equation}\label{eq:c}
\ct(\nu,J) = \Phi(\nu) + \sum_{\substack{1\le a\le n-1\\ 1\le i\le \mu_1}}
 |J_i^{(a)}|.
\end{equation}
Then the combinatorial analogue of \eqref{eq:sym-super} in terms
of rigged configurations is
\begin{equation*}
\St_{\la\mu}(q) = \sum_{(\nu,J)\in\RC(\la,\mu)} q^{\ct(\nu,J)}.
\end{equation*}

For the completely antisymmetric $q$-supernomials $\St'_{\la\mu}(q)$, set
$m_i^{(a)}=\nu_i^{(a)}-\nu_{i+1}^{(a+1)}$ and
$p_i^{(a)}=\nu_i^{(a+1)}-\nu_i^{(a)}$. In this case attach labels
to the last $m_i^{(a)}$ boxes in the $i$-th row of $\nu^{(a)}$.
A label $j$ in row $i$ of $\nu^{(a)}$ should satisfy $0\le j\le p_i^{(a)}$.
If $j=p_i^{(a)}$ the quantum number $j$ is called singular.
As in the symmetric case, riggings which differ by reordering of the labels 
within the same row of $\nu^{(a)}$ are identified. That is, to each 
$\nu^{(a)}$ and row $i$ one may associate a partition $J_i^{(a)}$ which 
lies in a box of width $m_i^{(a)}$ and height $p_i^{(a)}$.

\begin{example}
The configurations of example \ref{ex:config anti} admit the following 
riggings, for example,
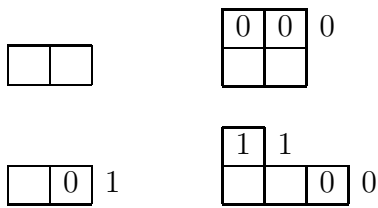
\begin{equation*}
\begin{array}{ll}
 \begin{array}[b]{|c|c|} \hline \spa & \spa \\ 
  \hline \end{array}
& \qquad
 \begin{array}[b]{|c|c|l} \cline{1-2} 0 & 0 & 0\\ 
  \cline{1-2} \spa & \spa & \\ \cline{1-2} \end{array}\\[4mm]
 \begin{array}[b]{|c|c|l} \cline{1-2} \spa & 0 & 1 \\ 
  \cline{1-2} \end{array}
& \qquad
 \begin{array}[b]{|c|c|c|l} 
  \cline{1-1} 1 & \multicolumn{3}{l}{1} \\
  \cline{1-3} \spa & \spa & 0 & 0 \\ \cline{1-3} \end{array}
\end{array}
\end{equation*}
where the vacancy numbers are denoted to the right of each row.
\end{example}
Denote the set of all antisymmetric rigged configurations
by $\RC'(\la,\mu)$. The statistic associated with 
$(\nu,J)\in \RC'(\la,\mu)$ is
\begin{equation}\label{eq:c'}
\ct'(\nu,J) = \sum_{\substack{1\le a\le n-1\\ 1\le i\le \mu_1}}
 |J_i^{(a)}|.
\end{equation}
Then the combinatorial analogue of \eqref{eq:anti-super} in terms
of rigged configurations is
\begin{equation*}
\St'_{\la\mu}(q) = \sum_{(\nu,J)\in\RC'(\la,\mu)} q^{\ct'(\nu,J)}.
\end{equation*}

\section{Ribbon tableaux}
\label{sec:RT}

Lascoux, Leclerc and Thibon \cite{LLT:1997} defined $q$-supernomial
coefficients in terms of ribbon tableaux. Ribbon tableaux are the
natural objects in the combinatorial description of the power-sum
plethysm operators on symmetric functions. Here we give a brief
review of the definition of the $q$-supernomial coefficients in terms 
of ribbon tableaux \cite{LLT:1997} and the reformulation in terms of 
inversion statistic on tuples of Young tableaux \cite{SSW:2000}.

\subsection{Cospin statistic}

To define $L$-ribbon tableaux we will mimic the definition of 
Young tableaux as chains in Young's lattice. We begin with the
review of the definition of Young tableaux.

Containment defines a partial order on partitions. If $\nu$ is
contained in $\la$ we write $\nu\subset \la$.
Young's lattice $\Y$ is the set of partitions under the partial order
$\subset$.
A horizontal strip is a skew shape that has at most one cell
in each column. A Young tableau of shape $\la/\nu$ and weight 
$\mu=(\mu_1,\ldots,\mu_r)$ is a chain of partitions
\begin{equation*}
\nu = \alpha^0 \subset \alpha^1 \subset \cdots \subset \alpha^r = \la
\end{equation*}
such that $\alpha^i/\alpha^{i-1}$ is a horizontal strip with $\mu_i$ cells.
As before, a Young tableau is represented graphically by the diagram of 
$\la/\nu$ where the cells in $\alpha^i/\alpha^{i-1}$ are numbered with $i$.
The set of all Young tableaux of shape $\la/\nu$ and weight $\mu$
is denoted by $\Tab(\la/\nu,\mu)$.

An $L$-ribbon is a connected skew shape consisting of $L$ cells, which 
does not contain any $2\times 2$ squares. The rightmost and lowermost cell 
is called the origin of the ribbon. Define the spin of a ribbon $R$ by 
$\spin(R)=\height(R)-1$ where $\height(R)$ is the height of $R$.

The relation $\nu\ribcov \mu$ on $\Y$ means that $\nu\subset \mu$
and the skew shape $\mu/\nu$ is an $L$-ribbon; it is the covering
relation for a partial order $\rible$ on $\Y$. Each component of
the poset $(\Y,\rible)$ has a unique minimum. The $\rible$-minima
are called $L$-cores. The $L$-ribbon lattice $\RH_L$ is the
component of the empty partition $\emptyset$ in $(\Y,\rible)$.

The skew shape $\mu/\nu$ is a horizontal $L$-ribbon strip of weight
$m$ if there is a saturated chain $\nu=\alpha^0\ribcov \alpha^1
\ribcov \cdots \ribcov \alpha^m=\mu$ such that the origin of each 
ribbon $R_i=\alpha^i/\alpha^{i-1}$ is in the lowermost cell 
of its column in $\mu/\nu$. If such a saturated chain exists it is unique if 
we require in addition that the origin of $R_i$ is to the
right of $R_{i-1}$.

An $L$-ribbon tableau $T$ of shape $\mu/\nu$ and weight 
$\la=(\la_1,\ldots,\la_n)$ is a chain of partitions
\begin{equation*}
\nu=\alpha^0 \rible \alpha^1 \rible \cdots \rible \alpha^n = \mu
\end{equation*}
such that $\alpha^i/\alpha^{i-1}$ is a horizontal $L$-ribbon strip
of weight $\la_i$. Graphically, $T$ may be represented by its
$L$-ribbons where each $L$-ribbon in $\alpha^i/\alpha^{i-1}$ is numbered
by $i$. The spin of an $L$-ribbon tableau $T$ is
the sum of the spins of its ribbons. Let $\Tab_L(\mu/\nu,\la)$ be
the set of all $L$-ribbon tableaux of shape $\mu/\nu$ and weight $\la$.
The cospin of $T\in \Tab_L(\mu/\nu,\la)$ is defined as
\begin{equation*}
\begin{split}
\maxspin(\mu/\nu)&=\max\{\spin(S) \mid S\in\Tab_L(\mu/\nu,\cdot)\}\\
\cospin(T)&=\frac{1}{2}(\maxspin(\mu/\nu)-\spin(T)).
\end{split}
\end{equation*}

\begin{example}\label{ex:3-ribbon}
The 3-ribbon tableau $T$ 
\begin{equation*}
\begin{array}{|c|c|c|c|c|c|}
\cline{1-3} &&4&\multicolumn{3}{c}{}\\
\cline{4-4} &3&&&\multicolumn{2}{c}{}\\
2&&&4&\multicolumn{2}{c}{}\\
\cline{1-3} &\multicolumn{1}{c}{}&2&&\multicolumn{2}{c}{}\\
\cline{2-2}\cline{4-6}&1&&&\multicolumn{1}{c}{3}&\\
\cline{3-3}\cline{5-5} 1&\multicolumn{1}{c}{}&&\multicolumn{1}{l}{3}&&\spa\\
\hline
\end{array}
\end{equation*}
of shape $\mu=(6,6,4,4,4,3)$ and weight $\la=(2,2,3,2)$
has $\spin(T)=14$ and $\cospin(T)=(16-14)/2=1$.
\end{example}

A standard $L$-ribbon tableau of shape $\mu/\nu$ is a saturated 
$\rible$-chain from $\nu$ to $\mu$. Equivalently, it is an $L$-ribbon
tableau of weight $(1^{|\mu/\nu|/L})$. The standardization $\st(T)$
of an $L$-ribbon tableau $T$ is the standard $L$-ribbon tableau obtained
by joining together the saturated chains of all its horizontal $L$-ribbon
strips.

\begin{example}
The standardization $\st(T)$ of the tableau in example \ref{ex:3-ribbon} is
\begin{equation*}
\begin{array}{|c|c|c|c|c|c|}
\cline{1-3} &&8&\multicolumn{3}{c}{}\\
\cline{4-4} &5&&&\multicolumn{2}{c}{}\\
3&&&9&\multicolumn{2}{c}{}\\
\cline{1-3} &\multicolumn{1}{c}{}&4&&\multicolumn{2}{c}{}\\
\cline{2-2}\cline{4-6}&2&&&\multicolumn{1}{c}{7}&\\
\cline{3-3}\cline{5-5} 1&\multicolumn{1}{c}{}&&\multicolumn{1}{l}{6}&&\spa\\
\hline
\end{array}
\end{equation*}
\end{example}

The generating function of $L$-ribbon tableaux with spin and cospin
statistics are defined as follows
\begin{equation}\label{eq:gen ribbon}
\begin{split}
\R_L(\mu/\nu,\la) &=\sum_{T\in \Tab_L(\mu/\nu,\la)} q^{\spin(T)}\\
\Rt_L(\mu/\nu,\la) &=\sum_{T\in \Tab_L(\mu/\nu,\la)} q^{\cospin(T)}.
\end{split}
\end{equation}

\subsection{Inversion statistic}
\label{sec:inv stat}

As was mentioned earlier, the set of $L$-ribbon tableaux can be
identified with the set of $L$-tuples of Young tableaux.
We follow mainly ref. \cite{SSW:2000} in this exposition.

Denote by $\Y^L$ the $L$-fold direct product of the poset $\Y$,
which is by definition the set of $L$-tuples of partitions
$\mud=(\mu^0,\mu^1,\dotsc,\mu^{L-1})$, also called
$L$-multipartitions. The set $\Y^L$ has a partial order $\nud\subset\mud$
given by $\nu^i\subset\mu^i$ for all $0\le i\le L-1$. 
Say that $\mud/\nud$ is a horizontal $L$-multistrip of weight $m$ if it is 
an $L$-tuple of horizontal strips with a total number of $m$ cells. 
Then an $L$-multitableau of shape $\mud/\nud$ and weight 
$\la=(\la_1,\ldots,\la_n)$ is a chain
\begin{equation}\label{eq:multi chain}
\nud = \alpha^{\bullet 0} \subset \alpha^{\bullet 1} \subset
 \cdots \subset \alpha^{\bullet n} = \mud
\end{equation}
such that $\alpha^{\bullet i}/\alpha^{\bullet i-1}$ is a horizontal
$L$-multistrip of weight $\la_i$. Equivalently, an $L$-multitableau $T$ 
may be viewed as an $L$-tuple of Young tableau $T=(T^0,T^1,\ldots,T^{L-1})$.
The set of all $L$-multitableaux of shape $\mud/\nud$ and weight $\la$
is denoted by $\Tab^L(\mud/\nud,\la)$.

Littlewood's $L$-quotient map gives a poset isomorphism
$\quot_L:\RH_L\rightarrow \Y^L$ (see \cite{FS:1998}). 
We use the following normalization of the $L$-quotient. Let $\la$ be a
partition. By appending zeros it may be assumed that the number of
parts of $\la$ is a multiple $m L$ of $L$. Define the staircase
partition $\rho^{(r)}=(r-1,r-2,\dots,1,0)\in\Z^r$. It can be shown
that $\la$ has empty $L$-core (that is, $\la\in\RH_L$) if and only
if $\la+\rho^{(mL)}$ has exactly $m$ parts that are congruent to
$i$ mod $L$ for each $0\le i\le L-1$. In that case, define the
$L$-quotient $\quot_L(\la)$ of $\la$ to be the $L$-tuple of
partitions $(\la^0,\la^1,\dots,\la^{L-1})$ where the
partition $L (\la^i+\rho^{(m)}) + i$ coincides with the
partition obtained by selecting the parts of $\la+\rho^{(mL)}$
that are congruent to $i$ modulo $L$. This definition is
independent of $m$.

The map $\quot_L$ induces a weight-preserving bijection
\begin{equation}\label{eq:SW}
\quot_L : \Tab_L(\mu,\la) \to \Tab^L(\mud,\la)
\end{equation}
from $L$-ribbon tableaux of shape $\mu$ to $L$-multitableaux of
shape $\mud=\quot_L(\mu)$.
This map is induced by observing that $\nu\rible \mu$ if and 
only if $\quot_L(\nu)\subset \quot_L(\mu)$, and $\la/\nu$
is a horizontal $L$-ribbon strip if and only if 
$\quot_L(\mu)/\quot_L(\nu)$ is a horizontal $L$-multistrip.
The bijection in \eqref{eq:SW} is called the Stanton--White 
bijection \cite{SW:1985}.

\begin{example}\label{ex:3-multi}
Under $\quot_3$ the 3-ribbon tableau $T$ of example \ref{ex:3-ribbon}
becomes
\begin{equation*}
\quot_3(T) = \left(\;
\begin{array}{|c|} \hline 2\\ \hline 1\\ \hline \end{array} \; , \;
\begin{array}{|c|c|} \hline 3&4\\ \hline 2&3\\ \hline \end{array} \; , \;
\begin{array}{|c|c|} \cline{1-1} 4&\multicolumn{1}{c}{}\\ \hline
 1&3\\ \hline \end{array}
\; \right).
\end{equation*}
\end{example}

A standard $L$-multitableau of shape $\mud/\nud$ is a saturated chain
in $\Y^L$ from $\nud$ to $\mud$. 
A cell in an $L$-multipartition $\mud$ is a triple
$s=(i,j,p)\in\Z^2\times\{0,1,\dots,L-1\}$ such that $(i,j)\in
\mu^p$. Write $\pos(s)$ for the index $p$ and set
$\row(s)=i$, $\col(s)=j$ and $\diag(s)=j-i$.
Define the partial order $\lecell$ on elements of
$\Z^2\times \{0,1,\dots,L-1\}$ by $s\lesscell s'$ if
$\diag(s)<\diag(s')$, or if $\diag(s)=\diag(s')$ and
$\pos(s)<\pos(s')$. 
The standardization of an $L$-multitableau $T$ of shape $\mud/\nud$
is the saturated chain in $\Y^L$ compatible with \eqref{eq:multi chain}
such that the cells within each horizontal $L$-multistrip are added
in $\lecell$-increasing order. 

\begin{example}\label{ex:stand}
The standardization of the 3-multitableau in example \ref{ex:3-multi} is
\begin{equation*}
\left(\;
\begin{array}{|c|} \hline 3\\ \hline 1\\ \hline \end{array} \; , \;
\begin{array}{|c|c|} \hline 5&9\\ \hline 4&6\\ \hline \end{array} \; , \;
\begin{array}{|c|c|} \cline{1-1} 8&\multicolumn{1}{c}{}\\ \hline
 2&7\\ \hline \end{array}
\; \right).
\end{equation*}
\end{example}

The set of $L$-ribbon tableaux is endowed with the spin and cospin
statistic. Under the bijection $\quot_L$ the cospin statistic
on $\Tab_L(\mu,\cdot)$ becomes an inversion statistic on
$\prod_{i=0}^{L-1} \Tab(\mu^i,\cdot)$ as shown in \cite{SSW:2000}.
We define the inversion statistic on standard $L$-multitableaux;
the inversion statistic of a general $L$-multitableau is the 
inversion of its standardization.

Let $T$ be a standard $L$-multitableau of shape $\mud$. For
$s\in\mud$ write $T(s)$ for the value of the cell $s$ in $T$,
which is the index $i$ such that $s$ is in the $i$-th multipartition 
in the chain $T$ but not in the $(i-1)$-st.

Given cells $s,t\in \mud$, say that $(s,t)$ is an
inversion of $T$ if the following conditions hold:
\begin{enumerate}
\item $\diag(s)=\diag(t)$ and $\pos(s)<\pos(t)$, or
$\diag(s)=\diag(t)-1$ and $\pos(s)>\pos(t)$. Note that in either
case $s\lesscell t$.
\item $\row(s) \le \row(t)$.
\item $T(t) < T(s) < T(\up{t})$ (where $\up{t}$ is the cell directly
above $t$ and $T(\up{t})=\infty$ if $\up{t}\not\in\mud$).
\end{enumerate}
Write $\inv(T)$ for the number of inversions of $T$.
A standard multitableau $T$ of shape $\mud=((1),\ldots,(1))$
can be identified with a permutation of the set $\{1,2,\dots,L\}$ 
and in this case $\inv(T)$ is the usual inversion number.

\begin{theorem}[\cite{SSW:2000}] \label{thm:cospin} 
For $T\in\Tab_L(\mu,\cdot)$,
\begin{equation*}
\cospin(T)=\inv(\quot_L(T)).
\end{equation*}
\end{theorem}

\begin{example}
The only inversion of the tableau $T$ of example \ref{ex:stand}
is $(s,t)$ where $T(s)=4$ and $T(t)=2$ so that $\inv(T)=1$ which agrees 
with the cospin of the 3-ribbon tableau of example \ref{ex:3-ribbon}.
\end{example}

\section{Statistic preserving bijection}
\label{sec:bij}

In this section we give a bijection from $L$-multitableaux
to rigged configurations in the symmetric (all tableaux in the
multitableaux are single rows) and the antisymmetric 
(all tableaux in the multitableaux are single columns) case.
This bijection preserves the statistics, which means that
the inversion statistic on multitableaux equals the statistic
on rigged configurations under this bijection.
In conjunction with the Stanton--White bijection of the previous
section, this defines a bijection from certain ribbon tableaux
to rigged configurations.

\subsection{Symmetric case}

Let $\mud=(\mu^0,\mu^1,\ldots,\mu^{L-1})$ be a multipartition
where the $\mu^i$ are single rows. The bijection
\begin{equation*}
\Psi: \Tab^L(\mud,\la) \to \RC(\la,\mu),
\end{equation*}
where $\mu$ is the partition with parts $|\mu^i|$, is defined recursively.
Let $T=(T^0,\ldots,T^{L-1})\in\Tab^L(\mud,\la)$ where
the entries of each tableau in this sequence are denoted
by $T^k=t^k_1t^k_2\ldots t^k_{|\mu^k|}$. The $L$-multitableau
$T$ is built up by successively adding the letters $t^k_i$
for $k=L-1,L-2,\ldots,0$ and $i=1,2,\ldots,|\mu^k|$.
The addition of the letter $t^k_i$ to $T$ corresponds to an
analogous operation on rigged configurations given 
by the following algorithm:
\begin{enumerate}
\item Define $(\nu,J)_{0,L-1}$ to be the empty rigged configuration
 in $\RC(\la,\mu)$.
\item Suppose the rigged configuration $(\nu,J)_{i-1,k}$ which corresponds 
 to the $L$-multitableau built up to letter $t^k_{i-1}$ is already known.
 In this notation we identify $(\nu,J)_{0,k}$ with $(\nu,J)_{|\mu^{k+1}|,k+1}$
 and set $t_0^k=1$.
 Denote by $(\nu,J)_{i-1,k}^{(a)}$ the $a$-th rigged partition in
 $(\nu,J)_{i-1,k}$. Adding the entry $t^k_i$ to the $k$-th tableau in $T$
 with $t^k_i\ge t^k_{i-1}$ corresponds to the following operation on the 
 rigged configuration $(\nu,J)_{i-1,k}$. For all $t^k_i\le a<n$, add a box to 
 row $i$ in $(\nu,J)_{i-1,k}^{(a)}$ such that a singular box is removed
 from row $i-1$ and the new box in row $i$ is singular.
 The resulting rigged configuration is $(\nu,J)_{i,k}$.
\item The rigged configuration $(\nu,J)=\Psi(T)$ is obtained from
 $(\nu,J)_{|\mu^0|,0}$ by inverting all quantum numbers. More precisely,
 a quantum number $j$ in row $i$ of the $a$-th rigged partition
 is replaced by $p_i^{(a)}-j$.
\end{enumerate}

\begin{figure}
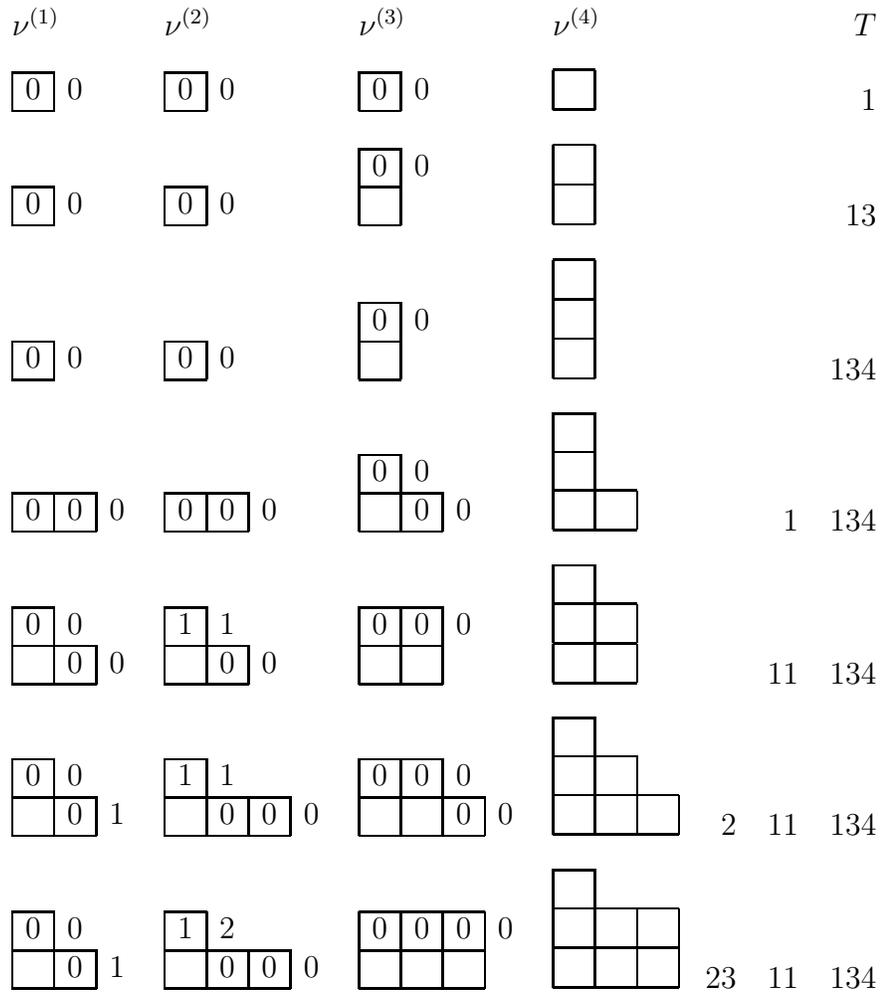

\begin{equation*}
\begin{array}{llllr}
\nu^{(1)}&\nu^{(2)}&\nu^{(3)}&\nu^{(4)}&T\\[3mm]
\begin{array}[b]{|c|l} \cline{1-1} 0&0 \\ \cline{1-1} \end{array}
&
\begin{array}[b]{|c|l} \cline{1-1} 0&0 \\ \cline{1-1} \end{array}
&
\begin{array}[b]{|c|l} \cline{1-1} 0&0 \\ \cline{1-1} \end{array}
&
\begin{array}[b]{|c|} \hline \phantom{0}\\ \hline \end{array}
&1\\[3mm]
\begin{array}[b]{|c|l} \cline{1-1} 0&0 \\ \cline{1-1} \end{array}
&
\begin{array}[b]{|c|l} \cline{1-1} 0&0 \\ \cline{1-1} \end{array}
&
\begin{array}[b]{|c|l} \cline{1-1} 0&0 \\ \cline{1-1} 
 &\\ \cline{1-1} \end{array}
&
\begin{array}[b]{|c|} \hline \phantom{0}\\ \hline \\ \hline \end{array}
&13\\[3mm]
\begin{array}[b]{|c|l} \cline{1-1} 0&0 \\ \cline{1-1} \end{array}
&
\begin{array}[b]{|c|l} \cline{1-1} 0&0 \\ \cline{1-1} \end{array}
&
\begin{array}[b]{|c|l} \cline{1-1} 0&0 \\ \cline{1-1} 
 &\\ \cline{1-1} \end{array}
&
\begin{array}[b]{|c|} \hline \phantom{0}\\ \hline \\ \hline \\
 \hline \end{array}
&134\\[3mm]
\begin{array}[b]{|c|c|l} \cline{1-2} 0&0&0\\ \cline{1-2} \end{array}
&
\begin{array}[b]{|c|c|l} \cline{1-2} 0&0&0\\ \cline{1-2} \end{array}
&
\begin{array}[b]{|c|c|l} \cline{1-1} 0&\multicolumn{2}{l}{0}\\
 \cline{1-2} &0&0\\ \cline{1-2} \end{array}
&
\begin{array}[b]{|c|c|} \cline{1-1}&\multicolumn{1}{c}{}\\
 \cline{1-1} &\multicolumn{1}{c}{}\\ \hline
 \phantom{1}&\phantom{1}\\ \hline \end{array}
&1\quad 134\\[3mm]
\begin{array}[b]{|c|c|l} \cline{1-1} 0&\multicolumn{2}{l}{0}\\
 \cline{1-2} &0&0\\ \cline{1-2} \end{array}
&
\begin{array}[b]{|c|c|l} \cline{1-1} 1&\multicolumn{2}{l}{1}\\
 \cline{1-2} &0&0\\ \cline{1-2} \end{array}
&
\begin{array}[b]{|c|c|l} \cline{1-2} 0&0&0\\ \cline{1-2}
 &&\\ \cline{1-2} \end{array}
&
\begin{array}[b]{|c|c|} \cline{1-1}&\multicolumn{1}{c}{}\\ \hline &\\ 
 \hline \phantom{1}&\phantom{1}\\ \hline \end{array}
&11\quad 134\\[3mm]
\begin{array}[b]{|c|c|l} \cline{1-1} 0&\multicolumn{2}{l}{0}\\
 \cline{1-2} &0&1\\ \cline{1-2} \end{array}
&
\begin{array}[b]{|c|c|c|l} \cline{1-1} 1&\multicolumn{3}{l}{1}\\
 \cline{1-3} &0&0&0\\ \cline{1-3} \end{array}
&
\begin{array}[b]{|c|c|c|l} \cline{1-2} 0&0&\multicolumn{2}{l}{0}\\ 
 \cline{1-3}&&0&0\\ \cline{1-3} \end{array}
&
\begin{array}[b]{|c|c|c|} \cline{1-1}&\multicolumn{2}{c}{}\\ 
 \cline{1-2} &&\multicolumn{1}{c}{}\\ 
 \hline \phantom{1}&\phantom{1}&\phantom{1}\\ \hline \end{array}
&2\quad 11\quad 134\\[3mm]
\begin{array}[b]{|c|c|l} \cline{1-1} 0&\multicolumn{2}{l}{0}\\
 \cline{1-2} &0&1\\ \cline{1-2} \end{array}
&
\begin{array}[b]{|c|c|c|l} \cline{1-1} 1&\multicolumn{3}{l}{2}\\
 \cline{1-3} &0&0&0\\ \cline{1-3} \end{array}
&
\begin{array}[b]{|c|c|c|l} \cline{1-3} 0&0&0&0\\
 \cline{1-3}&&&\\ \cline{1-3} \end{array}
&
\begin{array}[b]{|c|c|c|} \cline{1-1}&\multicolumn{2}{c}{}\\ 
 \hline &&\\
 \hline \phantom{1}&\phantom{1}&\phantom{1}\\ \hline \end{array}
&23\quad 11\quad 134
\end{array}
\end{equation*}
\caption{The intermediate steps in the algorithm for $\Psi(T)$ with $T$ of
example \ref{ex:rc sym}}
\label{figure:rc sym}
\end{figure}

\begin{example}\label{ex:rc sym}
Take the 3-multitableau 
$T=(\; \begin{array}{|c|c|} \hline 2&3\\ \hline \end{array} \;,\;
       \begin{array}{|c|c|} \hline 1&1\\ \hline \end{array} \;,\;
       \begin{array}{|c|c|c|} \hline 1&3&4\\ \hline \end{array}
   \;)$
which has inversion statistic $\inv(T)=3$. The corresponding 
rigged configuration $\Psi(T)$ is
\begin{equation*}
\begin{array}{|c|c|l} \cline{1-1} 0&\multicolumn{2}{l}{0}\\
 \cline{1-2} &1&1 \\ \cline{1-2} \end{array}
\qquad
\begin{array}{|c|c|c|l} \cline{1-1} 1&\multicolumn{3}{l}{2} \\
 \cline{1-3} &0&0&0 \\ \cline{1-3} \end{array}
\qquad
\begin{array}{|c|c|c|l} \cline{1-3} 0&0&0&0\\ \cline{1-3} &&&\\ 
 \cline{1-3} \end{array}
\end{equation*}
with $\ct(\Psi(T))=3$. All intermediate steps in the algorithm for the
calculation of $\Psi(T)$ are listed in Figure \ref{figure:rc sym}.
$\Psi(T)$ is obtained from the last rigged configuration by
inverting all quantum numbers.
\end{example}

\begin{theorem}
The map $\Psi:\Tab^L(\mud,\la) \to \RC(\la,\mu)$ defined by
the above algorithm is a statistic preserving bijection.
That is, $\inv(T)=\ct(\Psi(T))$ for $T\in\Tab^L(\mud,\la)$.
\end{theorem}

\begin{proof}
First we prove that $\Psi$ defines a bijection.
To this end it needs to be shown that $\Psi$ is well-defined
and that the algorithm is invertible.

Adding the letter $t$ to the $i$-th column of a tableau in $T$
adds a box in row $i$ to the $a$-th rigged partition for all
$t\le a<n$. This operation leaves all vacancy numbers unchanged,
except $p_i^{(t-1)}$ which changes by $+1$.
Since $t^k_i\ge t^k_{i-1}$ there hence exists a singular box in row
$i-1$ in $(\nu,J)_{i-1,k}^{(a)}$ for all $t_i^k\le a<n$.
This shows that $\Psi$ is well-defined.

It is easy to see that the algorithm for $\Psi$ is invertible.
Given $(\nu,J)\in\RC(\la,\mu)$ one constructs 
$T=(T^0,\ldots,T^{L-1})\in\Tab^L(\mud,\la)$ with 
$T^k=t_1^kt_2^k\cdots t_{|\mu^k|}^k$ by successively determining 
$t_i^k$ for $k=0,1,\ldots,L-1$ and $i=|\mu^k|,|\mu^k-1|,\ldots,1$.
Set $(\nu,J)_{|\mu^0|,0}$ to be 
$(\nu,J)$ with inverted quantum numbers and let $(\nu,J)_{i+1,k}$ be 
the rigged configuration after determining $t_{i+1}^k$.
To determine $t_i^k$, one finds the smallest index $t$ such that
there is a singular string in row $i$ of $(\nu,J)_{i+1,k}^{(a)}$ for all
$t\le a<n$. Set $t_i^k=t$. The rigged configuration $(\nu,J)_{i,k}$ is 
obtained from $(\nu,J)_{i+1,k}$ by removing the selected singular
strings in row $i$ and making the new boxes in row $i-1$ singular.
Then it is clear that $t_i^k\le t_{i+1}^k$
(where $t_{i+1}^k=\infty$ if $i=|\mu^k|$) so that $T$ will indeed
be an $L$-multipartition of shape $\mud$ and content $\la$.

It remains to prove that $\Psi$ preserves the statistic.
Since we are dealing with $L$-multitableaux of single
rows the inversion statistic of Section \ref{sec:inv stat}
simplifies as follows. The pair $(s,t)$ is an inversion of $T$
if $T(t)<T(s)$ and either $\diag(s)=\diag(t)$ and $\pos(s)<\pos(t)$
or $\diag(s)=\diag(t)-1$ and $\pos(s)>\pos(t)$.
Furthermore, the following relation between a multitableau $T$ and 
the corresponding configuration $\nu$ will be needed
\begin{equation}\label{eq:rel Tn}
\nu_i^{(a)}-\nu_i^{(a-1)}=\card\{s\in\mud \mid 
 \text{$T(s)=a$, $\diag(s)=i-1$}\}.
\end{equation}

Building $T$ up successively, the addition of the letter
$t^k_i$ in box $t\in\mud$ changes the inversion statistic by
\begin{equation*}
\begin{split}
\Delta\inv &= 
 \card \{s\in \mud \mid \text{$\pos(s)>\pos(t)$, $T(s)<T(t)$, 
 $\diag(s)=\diag(t)$}\}\\
 &+ \card \{s\in \mud \mid \text{$\pos(s)>\pos(t)$, $T(s)>T(t)$,
 $\diag(s)=\diag(t)-1$}\}.
\end{split}
\end{equation*}
Denote by $\nu$ the configuration corresponding to $T$ built up to 
letter $t_i^k$ in the algorithm for $\Psi$. 
Then, setting $a=T(t)$ and inserting \eqref{eq:rel Tn} yields
\begin{equation}\label{eq:c inv}
\begin{split}
\Delta\inv &= \sum_{\alpha=1}^{a-1} (\nu_i^{(\alpha)}-\nu_i^{(\alpha-1)})
 + \sum_{\alpha=a+1}^n (\nu_{i-1}^{(\alpha)}-\nu_{i-1}^{(\alpha-1)})\\ 
 &= \nu_i^{(a-1)}+\nu_{i-1}^{(n)}-\nu_{i-1}^{(a)}.
\end{split}
\end{equation}

On the other hand, the change of the rigged configuration
statistic $\ct$, taking into account that we have to invert
the quantum numbers, is given by
\begin{equation}\label{eq:c ct}
\begin{split}
\Delta \ct &= \nu_{i+1}^{(a-1)}+(\nu_{i-1}^{(a+1)}-\nu_{i-1}^{(a)})
 +(\nu_{i-1}^{(a+2)}-\nu_{i-1}^{(a+1)})+\cdots
 +(\nu_{i-1}^{(n)}-\nu_{i-1}^{(n-1)})\\
 &+(\nu_i^{(a-1)}-\nu_{i+1}^{(a-1)})\\
&=\nu_i^{(a-1)}+\nu_{i-1}^{(n)}-\nu_{i-1}^{(a)},
\end{split}
\end{equation}
where the first line comes from the change in $\Phi(\nu)$
and the second line is induced by the change of the vacancy number.
Since \eqref{eq:c inv} and \eqref{eq:c ct} agree at each step
in the algorithm and $\inv(\emptyset)=\ct(\emptyset)=0$, $\Psi$ is 
statistic preserving.
\end{proof}

\subsection{Antisymmetric case}

Let $\mud=(\mu^0,\mu^1,\ldots,\mu^{L-1})$ be a multipartition
where the $\mu^i$ are single columns. As in the symmetric case,
the bijection
\begin{equation*}
\Psi': \Tab^L(\mud,\la) \to \RC'(\la,\mu),
\end{equation*}
is defined recursively where $\mu$ is the partition with parts
$|\mu^i|$. Let $T=(T^0,\ldots,T^{L-1})\in\Tab^L(\mud,\la)$.
The letters in the single column tableau $T^k$ are denoted by 
$t^k_{|\mu^k|} t^k_{|\mu^k|-1} \cdots t^k_1$ in strictly decreasing order.
The $L$-multitableau $T$ is built up by successively adding the letters 
$t^k_i$ for $k=L-1,L-2,\ldots,0$ and $i=1,2,\ldots,|\mu^k|$.
The addition of the letter $t^k_i$ to $T$ corresponds to an
analogous operation on rigged configurations given 
by the following algorithm:
\begin{enumerate}
\item Define $(\nu,J)_{0,L-1}$ to be the empty rigged configuration
 in $\RC'(\la,\mu)$.
\item Suppose the rigged configuration $(\nu,J)_{i-1,k}$ which corresponds 
 to the $L$-multitableau built up to letter $t^k_{i-1}$ is already known.
 Here we identify $(\nu,J)_{0,k}$ with $(\nu,J)_{|\mu^{k+1}|,k+1}$ and
 set $t^k_0=0$.
 Denote by $(\nu,J)_{i-1,k}^{(a)}$ the $a$-th rigged partition in
 $(\nu,J)_{i-1,k}$. Adding the entry $t^k_i$ to the $k$-th tableau in $T$
 with $t^k_i>t^k_{i-1}$ corresponds to the following operation on the 
 rigged configuration $(\nu,J)_{i-1,k}$. For all $t^k_i\le a<n$, add a box to 
 row $i$ in $(\nu,J)_{i-1,k}^{(a)}$ and make the new box singular.
 Remove a singular quantum number from row $i-1$ in all $(\nu,J)_{i-1,k}^{(a)}$
 for $t_i^k-1\le a<n$. Call the resulting rigged configuration 
 $(\nu,J)_{i,k}$.
\item The rigged configuration $(\nu,J)=\Psi'(T)$ is obtained from
 $(\nu,J)_{|\mu^0|,0}$ by inverting all quantum numbers. More precisely,
 a quantum number $j$ in row $i$ of the $a$-th rigged partition
 is replaced by $p_i^{(a)}-j$.
\end{enumerate}

\begin{figure}
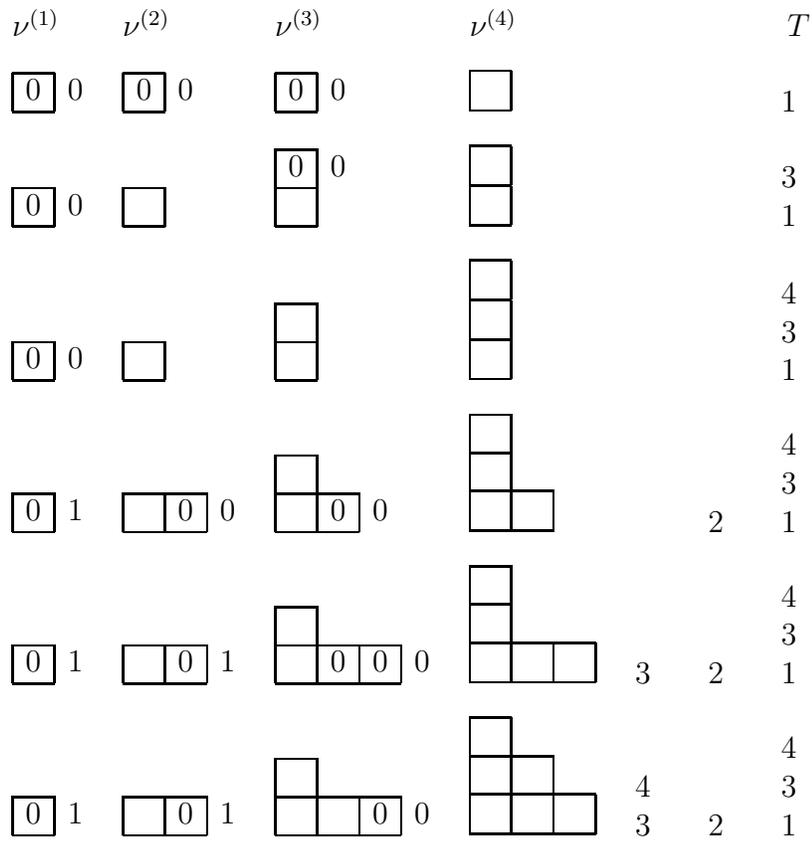

\begin{equation*}
\begin{array}{llllr}
\nu^{(1)}&\nu^{(2)}&\nu^{(3)}&\nu^{(4)}&T\\[3mm]
\begin{array}[b]{|c|l} \cline{1-1} 0&0 \\ \cline{1-1} \end{array}
&
\begin{array}[b]{|c|l} \cline{1-1} 0&0 \\ \cline{1-1} \end{array}
&
\begin{array}[b]{|c|l} \cline{1-1} 0&0 \\ \cline{1-1} \end{array}
&
\begin{array}[b]{|c|} \hline \phantom{0}\\ \hline \end{array}
&
\begin{array}[b]{c}1\end{array}\\[3mm]
\begin{array}[b]{|c|l} \cline{1-1} 0&0 \\ \cline{1-1} \end{array}
&
\begin{array}[b]{|c|l} \cline{1-1} \phantom{0}& \\ \cline{1-1} \end{array}
&
\begin{array}[b]{|c|l} \cline{1-1} 0&0 \\ \cline{1-1} &\\ 
  \cline{1-1} \end{array}
&
\begin{array}[b]{|c|} \hline \phantom{0}\\ \hline \\ \hline \end{array}
&
\begin{array}[b]{c} 3\\1 \end{array}\\[3mm]
\begin{array}[b]{|c|l} \cline{1-1} 0&0 \\ \cline{1-1} \end{array}
&
\begin{array}[b]{|c|l} \cline{1-1} \phantom{0}& \\ \cline{1-1} \end{array}
&
\begin{array}[b]{|c|l} \cline{1-1} \phantom{0}& \\ \cline{1-1} &\\ 
  \cline{1-1} \end{array}
&
\begin{array}[b]{|c|} \hline \phantom{0}\\ \hline \\ \hline \\
 \hline \end{array}
&
\begin{array}[b]{c}4\\3\\1\end{array}\\[3mm]
\begin{array}[b]{|c|l} \cline{1-1} 0&1 \\ \cline{1-1} \end{array}
&
\begin{array}[b]{|c|c|l} \cline{1-2} \phantom{0}&0&0\\ \cline{1-2} \end{array}
&
\begin{array}[b]{|c|c|l} \cline{1-1} \phantom{0}&\multicolumn{2}{l}{}\\
 \cline{1-2} &0&0\\ \cline{1-2} \end{array}
&
\begin{array}[b]{|c|c|} \cline{1-1}&\multicolumn{1}{c}{}\\
 \cline{1-1} &\multicolumn{1}{c}{}\\ \hline
 \phantom{1}&\phantom{1}\\ \hline \end{array}
&
\begin{array}[b]{c}2\end{array}\quad 
 \begin{array}[b]{c}4\\3\\1\end{array}\\[3mm]
\begin{array}[b]{|c|l} \cline{1-1} 0&1 \\ \cline{1-1} \end{array}
&
\begin{array}[b]{|c|c|l} \cline{1-2} \phantom{0}&0&1\\ \cline{1-2} \end{array}
&
\begin{array}[b]{|c|c|c|l} \cline{1-1} \phantom{0}&\multicolumn{3}{l}{}\\
 \cline{1-3} &0&0&0\\ \cline{1-3} \end{array}
&
\begin{array}[b]{|c|c|c|} \cline{1-1}&\multicolumn{2}{c}{}\\ 
 \cline{1-1} &\multicolumn{2}{c}{}\\ 
 \hline \phantom{1}&\phantom{1}&\phantom{1}\\ \hline \end{array}
&
\begin{array}[b]{c}3\end{array}\quad
 \begin{array}[b]{c}2\end{array}\quad
 \begin{array}[b]{c}4\\3\\1\end{array}\\[3mm]
\begin{array}[b]{|c|l} \cline{1-1} 0&1 \\ \cline{1-1} \end{array}
&
\begin{array}[b]{|c|c|l} \cline{1-2} \phantom{0}&0&1\\ \cline{1-2} \end{array}
&
\begin{array}[b]{|c|c|c|l} \cline{1-1} \phantom{0}&\multicolumn{3}{l}{}\\
 \cline{1-3} &\phantom{0}&0&0\\ \cline{1-3} \end{array}
&
\begin{array}[b]{|c|c|c|} \cline{1-1}&\multicolumn{2}{c}{}\\ 
 \cline{1-2} &&\multicolumn{1}{c}{}\\ 
 \hline \phantom{1}&\phantom{1}&\phantom{1}\\ \hline \end{array}
&
\begin{array}[b]{c}4\\3\end{array}\quad
 \begin{array}[b]{c}2\end{array}\quad
 \begin{array}[b]{c}4\\3\\1\end{array}
\end{array}
\end{equation*}
\caption{The intermediate steps in the algorithm for $\Psi'(T)$ with $T$ of
example \ref{ex:rc anti}}
\label{figure:rc anti}
\end{figure}

\begin{example}\label{ex:rc anti}
Consider the 3-multitableau
\begin{equation*}
T= \bigg( \;
\raisebox{-0.4cm}{
\begin{array}[b]{|c|} \hline 4\\ \hline 3\\ \hline \end{array} \;,\;
\begin{array}[b]{|c|} \hline 2\\ \hline \end{array} \;,\;
\begin{array}[b]{|c|} \hline 4\\ \hline 3\\ \hline 1\\ \hline \end{array}}
 \;\, \bigg)
\end{equation*}
with inversion statistic $\inv(T)=2$.
The corresponding rigged configuration $\Psi'(T)$ is
\begin{equation*}
\begin{array}[b]{|c|l} \cline{1-1} 1&1\\ \cline{1-1} \end{array}
\qquad
\begin{array}[b]{|c|c|l} \cline{1-2} \phantom{1}&1&1\\ \cline{1-2} \end{array}
\qquad
\begin{array}[b]{|c|c|c|l} \cline{1-1} &\multicolumn{3}{c}{}\\
 \cline{1-3} \phantom{1}&\phantom{1}&0&0\\ \cline{1-3} \end{array}
\end{equation*}
with $\ct'(\Psi'(T))=2$. The intermediate steps of the algorithm
are listed in Figure \ref{figure:rc anti} and $\Psi'(T)$ is obtained
from the last line by inversion of the quantum numbers.
\end{example}

\begin{theorem}
The map $\Psi':\Tab^L(\mud,\la) \to \RC'(\la,\mu)$ defined by
the above algorithm is a statistic preserving bijection.
That is, $\inv(T)=\ct'(\Psi'(T))$ for $T\in\Tab^L(\mud,\la)$.
\end{theorem}

\begin{proof}
We begin by proving that $\Psi'$ defines a bijection.
Adding the letter $t$ to the $i$-th row of a tableau in $T$
adds a box in row $i$ to the $a$-th rigged partition for all
$t\le a<n$. This operation leaves all vacancy numbers unchanged,
except $p_i^{(t-1)}$ which changes by $+1$.
Since $t^k_i> t^k_{i-1}$ there exists a singular box in row
$i-1$ in $(\nu,J)_{i-1,k}^{(a)}$ for all $t_i^k-1\le a<n$.
This shows that $\Psi'$ is well-defined.
It is easy to see that the algorithm for $\Psi'$ is invertible.

It remains to show that the statistic is preserved.
In the single column case $(s,t)$ is an inversion if
$\diag(s)=\diag(t)$, $\pos(s)<\pos(t)$ and $T(t)<T(s)<T(\up{t})$.
Furthermore, the analog of \eqref{eq:rel Tn} becomes
\begin{equation}\label{eq:rel Tn anti}
\nu^{(a)}_i-\nu^{(a-1)}_i = \card\{s\in\mud \mid T(s)=a, \diag(s)=1-i\}.
\end{equation}

Adding the letter $t_i^k$ to $T$ in box $t\in\mud$ changes the inversion
number by
\begin{equation*}
\begin{split}
\Delta \inv &= \card\{ s\in\mud \mid \text{$\pos(s)>\pos(t)$,
 $T(s)<T(t)$, $\diag(s)=\diag(t)$}\}\\
 &- \card\{ s\in\mud \mid \text{$\pos(s)>\pos(t)$, $T(s)\le T(t)$,
 $\diag(s)=\diag(t)-1$}\}
\end{split}
\end{equation*}
where the last line subtracts all contributions which do not
satisfy $T(t)<T(\up{s})$. Denote by $\nu$ the configuration
corresponding to $T$ built up to letter $t_i^k$ in the
algorithm $\Psi'$. Then, setting $a=T(t)$ and inserting
\eqref{eq:rel Tn anti} yields
\begin{equation*}
\begin{split}
\Delta \inv &= \sum_{\alpha=1}^{a-1}(\nu_i^{(\alpha)}-\nu_i^{(\alpha-1)})
 - \sum_{\alpha=1}^a (\nu_{i+1}^{(\alpha)}-\nu_{i+1}^{(\alpha-1)})\\
 &= \nu_i^{(a-1)}-\nu_{i+1}^{(a)}.
\end{split}
\end{equation*}
On the other hand, taking into account that the quantum numbers need
to be inverted for the correct statistic, we have
\begin{equation*}
\Delta \ct' = \nu_i^{(a-1)} - \nu_{i+1}^{(a)}
\end{equation*}
from the change in vacancy numbers. Since $\Delta \inv$ and
$\Delta \ct'$ agree at each step of the algorithm, $\Psi'$ is
statistic preserving.
\end{proof}

\subsection*{Acknowledgments}
I would like to thank the Max-Planck-Institut f\"ur Mathematik
in Bonn for hospitality where this work was completed.

\bibliographystyle{plain}

\vskip 2pc
\begin{tabular}{ll}
 Department of Mathematics & Max-Planck-Institut\\
 University of California  & f\"ur Mathematik in Bonn\\
 One Shields Ave           & Vivatsgasse 7\\
 Davis, CA 95616-8633      & 53111 Bonn\\
 U.S.A.                    & Germany\\
 e-mail: {\tt anne@math.ucdavis.edu} & e-mail: {\tt annes@mpim-bonn.mpg.de}
\end{tabular}
\end{document}